[1994/12/01]
\documentclass[reqno,13pt]{amsart}
\usepackage{bbm}
\def\1{\mathbbm{1}}

\usepackage{amsmath,amsthm}
\usepackage[mathscr]{euscript}
\usepackage{mathrsfs}
\usepackage{chngcntr}
\usepackage{graphicx}
\usepackage{amsmath,amssymb}
\usepackage{circuitikz}

\def\A{{\mathbb{A}}}
\def\R{{\mathbb{R}}}

\def\al{\alpha}
\def\om{\omega}

\def\ga{\gamma}
\def\si{\sigma}

\def\la{\lambda}

\def\t{\tau}

\def\B{\mathbb{B}}

\def\calL{{\mathcal{L}}}

\def\calA{{\mathcal{A}}}

\def\R{\mathbb R}

\def\C{\mathbb C}

\def\N{\mathbb N}

\def\M{\mathbb M}
\def\F{\mathbb F}

\def\L{\mathbb L}

\def\G{\mathbb G}

\newtheorem{definition}{Definition}[section]

\newtheorem{lemma}{Lemma}[section]
\newtheorem{proposition}{Proposition}[section]
\newtheorem{corollary}{Corollary}[section]
\newtheorem{theorem}{Theorem}[section]
\theoremstyle{remark}
\newtheorem{remark}{Remark}[section]
\newtheorem{example}{Example}[section]
\newtheorem{mainassumptions}{Main Assumptions}[section]

\title{Controllability of vertex delay type problems by the regular linear systems approach}
\author{ Y. El Gantouh, S. Hadd and A. Rhandi}
\address{ Yassine El gantouh, Departamento de Matem\'{a}ticas, Universidad Aut\'{o}noma de Madrid, 28049 Madrid, Spain; elgantouhyassine@gmail.com}
\address{Said Hadd: Department of Mathematics, Faculty of Sciences, Ibn Zohr University,  Agadir, Morocco; s.hadd@uiz.ac.ma}
\address{Rhandi: Dipartimento di Ingegneria dell'Informazione, Ingegneria Elettrica e Matematica Applicata, Università degli Studi di Salerno, Via Giovanni Paolo II, 132, 84084 Fisciano (Sa), Italy}
\thanks{This work has been supported by COST Action CA18232.}
\thanks{}
\keywords{network systems, delay systems, boundary control, Kalman-type conditions, structural controllability.}

\medskip

\begin{document}
	\maketitle
	
	\renewcommand{\sectionmark}[1]{}
	\begin{abstract}
In this paper, we study the well-posedness and approximate controllability of a class of network systems having delays and controls at the boundary conditions. The particularity of this work is that the network system is defined on infinite metric graphs. This fact offers many difficulties in applying the usual methods. In fact, the well-posedness of the delay network system is obtained by using a semigroup approach on product spaces which is based on the concept of feedback theory of infinite-dimensional linear systems. This technique allows us to reformulate the delay system into a free-delay distributed control system. From this transformation,  we deduce necessary and sufficient conditions for the boundary approximate controllability of such systems. Furthermore, a Rank condition to the approximate controllability is also obtained. This condition coincides with the usual Kalman controllability criterion in the case of a simple transport process on a finite graph. Finally, by applying our approach to a linear Eulerian model (with airborne delays) for (ATFM), we provide a new algebraic condition for the controllability of such a model in terms of the generic rank of the so-called extended controllability matrix.
	\end{abstract}

	\section{Introduction}
\label{sec:introduction}
Consider the following system of network transport equations
\begin{align*}
\mathsf{(vdp)}\begin{cases}
\dfrac{\partial }{\partial t}z_{j}(t,x)= c_{j}(x)\dfrac{\partial }{\partial x}z_{j}(t,x)+q_{j}(x).z_{j}(t,x),&x\in (0,1), \; t\geq 0, \\  
z_{j}(0,x)= g_{j}(x), &x\in (0,1),\quad(\text{IC}) \\
\mathsf{i}^{-}_{ij}c_{j}(1)z_{j}(t,1)= \mathsf{w}_{ij}\displaystyle\sum_{k\in \mathbb{N}} \mathsf{i}^{+}_{ik}\big[c_{k}(0)z_{k}(t,0)+ L_{k}z_{k}(\cdot +t,\cdot)\big]+\mathsf{b}_{il}u_{l}(t),& t\geq 0,\qquad\quad (\text{BC})\\
z_j(\theta,x)=\varphi_{j}(\theta,x) ,&x\in (0,1),\; \theta\in[-r,0],          
\end{cases}
\end{align*}
for $ i,\,j\in \mathbb{N}$ and $l\in \{1,\ldots,N\}$. This system  is a macroscopic model describing the continuous evolution of flows in a network subjected to transmission conditions with time delayed behavior. The corresponding transport equations are defined on the edges of an infinite metric graph $\mathsf{G}$, whose edges are identified with a collection of intervals with endpoints "glued" to the graph structure. The connection of such edges being described by the coefficients $\textsf{i}^{-}_{ij},\textsf{i}^{+}_{ik}\in \{0,1\}$ for $i,j,k\in\mathbb{N}$. The flow velocity along an edge $ e_j $ is determined by the function $ c_j $, whereas its absorption is determined by the function $ q_j $. The boundary condition (BC) determines the propagation of flows along the various components of the network, where the weights $w_{ij}\in \mathbb{R}_+$ and the vertex delay operators $L_k$ (see \eqref{S1.3} below), for $i,j,k\in\mathbb{N} $, express the proportion of mass being redistributed into the edges and the hereditary effects of the transmission conditions at the vertices (exteriors), respectively. Moreover, for $i,l\in\mathbb{N}\times \{1,\ldots,N\}$, the coefficients $\mathsf{b}_{il}$ denotes the entries of the so-called input matrix $ K $ and $ u_l $ define the control functions at the vertices, whose acts by adjusting the distribution of materials throughout the edges of the underlying network. 

In recent years, advanced tools on control theory have been applied to the analysis and control of PDEs on networks. (see e.g., \cite{DZ,LLS,Zua,HLS,EKNS,EKKNS,Br,EHR} and references therein). The goal was to use a sophisticated mathematical background to address well-posedness and control properties of transport equations such as controllability and observability \cite{DZ}, stabilization \cite{Zua}, switching controllability \cite{HLS}, and boundary controllability \cite{EKNS,EHR,EH}. In the present work, we are interested in studying the well-posedness and the approximate controllability of a transport
tree-like network with infinitely many edges including hereditary effects in the transmission conditions of the form $\mathsf{(vdp)}$.

In the absence of delays in the transmission conditions (i.e $L_k\equiv 0$), the  well-posedness of $\mathsf{(vdp)}$ is studied in \cite{BDK,Dorn-infinie,EKNS,EHR, KMN,TE}, where the corresponding graph is assumed to be finite. The technique used in these papers is based on Greiner's approach (see \cite{Gr}) which deals with equation with perturbed boundary conditions. In particular, the work \cite{EHR} uses a perturbation theorem developed in \cite{HMR} to prove the well-posedness of the network system. 

In the presence of a delay term  at the boundary conditions and the absence of control forces ($\mathsf{b}_{il}=0$), the authors of \cite{BDR12} proved the well-posedness and stability property of the equation $\mathsf{(vdp)}$ in the case of a particular bounded delay operator and a finite graph. However, when the delay is distributed, the problem of the existence and uniqueness of solutions of $\mathsf{(vdp)}$ presents some difficulties. Indeed, it is not clear how to use the classical semigroup theorems such as Hille-Yosida  and/or Lumer-Phillips theorems (see e.g. \cite{EN}) to prove the existence of solutions for the problem $\mathsf{(vdp)}$. The first main contribution of the present work is to rely on the paper \cite{HMR} to introduce a semigroup approach on product spaces to reformulate the delay network system $\mathsf{(vdp)}$ as an infinite dimensional well-posed distributed linear open-loop systems in the sense of \cite{Sa,St,WC}.  In fact, by introducing suitable conditions on the coefficients of the systems $\mathsf{(vdp)}$, see Section \ref{S.3}, we reformulate the delay system  $\mathsf{(vdp)}$ as the system \eqref{open-loop} (see Theorem \ref{Main-section3}). This will bring the controllability of $\mathsf{(vdp)}$ in line with the standard controllability definitions in systems theory, see e.g. \cite{TW}.  We mention that for finite-dimensional linear systems,  exact controllability and approximate controllability coincide, and are characterized by the famous Kalman rank condition and the condition of Hautus \cite{TSH}.

After having established the well-posed character of the system $\mathsf{(vdp)}$,  in the second part of this article (see Section \ref{S.4}) we will study the concept of (infinite-time) controllability for a such system. It should be noted that the system $\mathsf{(vdp)}$ is an infinite dimensional transport equation, so  approximate and exact controllability are not the same for the this system. Usually, when a system is affected by delays, it is more practical to study approximate controllability. Thus, necessary and sufficient conditions for the approximate controllability of the delay network system $\mathsf{(vdp)}$ are introduced using a duality and Laplace transform arguments. In particular we propose a controllability criterion in terms of a Kalman-type rank condition for the delay system $\mathsf{(vdp)}$  involving the graph structure. Let us mention that, our approach can not be applied to treat the controllability in finite time. However, as $\mathsf{(vdp)}$ is equivalent to a distributed linear system governed by a semigroup and an admissible control operator, one can then find conditions for which the system is approximately controllable in finite time by using Gramian's concept of observability, the fact that controllability and observability are dual properties, and a perturbation argument.  

We recall that the parameters of a  network system  are not precisely known (independent free parameters), see e.g. \cite{LSB,LBa}. The controllability corresponding to this case is called the  \emph{structural controllability}, see e.g. \cite{LBa,Li,LSB,Ol,SP,PKA,LiA} for the definitions and properties.

One of the main problems that arise in the study of the controllability of transport network systems is the choice of an appropriate boundary control abstract framework. The latter depends on the appropriate choice of transmission conditions between the various components of the underlying network. For example, if one chooses the transmission conditions to present the standard conditions \emph{Kirchhoff}, this can be seen as imposing additional constraints on the flows. This fact has been shown to have an effect on the propagation of streams along the network, see for example \cite{Br,EKNS}. The authors of \cite{EKNS} observed that only a certain subset of mass distributions can be obtained on the edges. This observation allows them to introduce a new concept of controllability called \emph{maximum controllability}. In fact, it has been proven that the flow can be maximally controlled at every vertex, at some vertices only, or at none of the vertices in the network. However, in the works cited above, only particular situations have been studied. Moreover, no condition on the structure of the graph nor the description of the flow (in terms of addition or subtraction of matter in the controlled vertex) was obtained.

It is remarkable that few works in the literature are devoted to characterizing the controllability properties of transport network systems. It seems that the attention was mainly directed to the long-term behavior and also to provide more accurate physical models describing the relationship between the graph topology, the dynamics and the boundary conditions imposed on the vertices of the networks (see, for example, \cite{BaN,JMZ,KMN} and their references).

We summarize the results of this paper as follows: Section \ref{S.2} is devoted to recalling the concept of feedback theory of infinite-dimensional regular linear systems as well as a key perturbation theorem (see Theorem \ref{nonhomo-bound}). Section \ref{S.3} is divided into two subsections. In fact, in Subsection \ref{sec3-sub1}, we introduce the main required assumptions, some notations on spaces and operators, and then use a product state space to reformulate delay network systems $\mathsf{(vdp)}$ into a perturbed boundary control system without delay \eqref{Sy.6}. In Subsection \ref{sec3-sub2}, we prove two main results, the first concerns a generation theorem (see Theorem \ref{T.2}), and the second concerns the transformation of the system \eqref{Sy.6} as an open-loop system \eqref{open-loop} (see Theorem \ref{Main-section3}), which implies the well-posedness of $\mathsf{(vdp)}$. In Section \ref{S.4}, based on theorem \ref{Main-section3}, we first define three types of approximate controllability for the network system $\mathsf{(vdp)}$ (see Definition \ref{types-controllability}), then we characterize each of these types (see Theorem \ref{Main-controllanility1} and Proposition \ref{partial-prop}). Moreover, in the Remark \ref{relation-between-controllability} we show the link that exists between these types of controllability. Furthermore, we will prove a rank condition for the approximate controllability of the system $\mathsf{(vdp)}$ (see Theorem \ref{T.4}). Section \ref{S.5} is about an application to an  Air Traffic Flow Management.

\subsection{Notation and Terminology}
Throughout the paper, $\mathbb{C},\mathbb{R},\mathbb{N},\mathbb{Q}$ are sets of complex, real, natural, and rational numbers, respectively. The cardinality of a finite set $\mathfrak{N}$ is denoted by $\#\mathfrak{N}$. $\ell^1:= \{(y_k)_{k\in \mathbb{N}}:\; \sum_{k=1}^{\infty} \vert y_k\vert<\infty \}$ is the space of all absolutely summable sequences of real or complex numbers. For an infinite matrix $D$, $(D)_{ij}$ indicates the element of $D$ which is located at its $i$-th row and $j$-th column, $D^{\top}$ is used to denote its matrix transpose, and Rg$\,(D)$ is used to denotes its range. If $D$ is an infinite matrix and $y=(y_j)$ is an infinite sequence (i.e., vector), we define $Dy$ by $(Dy)_i:=\sum_{j=1}^{\infty}(D)_{ij}y_j$, for each $i\in \mathbb{N}$, for which this infinite series converges.

Let $\mathsf{G}$ be an infinite connected metric graph in the sense that $ (\mathsf{G},d) $ is a connected metric space for which there exists a countable set $\mathsf{V}:=\{\mathsf{v}_i, \, i\in \mathsf{I} \}$, the set of vertices, and a partition $ \{\mathsf{e}_j, \, j\in \mathsf{J} \} $ of $ \mathsf{G}\setminus \mathsf{V} $ with $ \mathsf{J} $ is a countable set (i.e $ \mathsf{G}\setminus \mathsf{V}= \cup_{j\in\mathsf{J}}\mathsf{e}_j $ and for $ j\neq j' $, $ \mathsf{e}_j \cap\mathsf{e}_{j'}=\emptyset $) such that for all $ j\in\mathsf{J},\, \mathsf{e}_j $ is isometric to an interval $ (0,l_j] $ with $ l_j< +\infty $. We call $\mathsf{e}_j $ an edge, $ l_j$ its length and denote by $ \mathsf{E}:=\{\mathsf{e}_j, \, j\in \mathsf{J} \} $ the set of all edges on $\mathsf{G}$, see e.g., \cite{Delio}. Note that if an edge $\mathsf{e}\in \mathsf{E}$ connects $\mathsf{v}_i,\mathsf{v}_j\in \mathsf{V}$, we say that $\mathsf{v}_i,\mathsf{v}_j$ are incident with $\mathsf{e}$ and we write $\mathsf{v}_i \sim \mathsf{v}_j$.
We use $\mathsf{e}=(\mathsf{v}_i,\mathsf{v}_j)$ to denote direct edges of $\mathsf{E}$, where the vertices $\mathsf{v}_i,\mathsf{v}_j$ represent its tail and head, respectively. In this case, $\mathsf{e}$ is called an outgoing edge of the vertex $\mathsf{v}_i$, whereas it is an incoming edge of vertex $\mathsf{v}_j$. The edge $ (\mathsf{v}_i,\mathsf{v}_i) $ is an example of a loop, and the vertex $\mathsf{v}$ that has no incident edges is called an isolated vertex. When a graph contains no loops it is called \emph{loop-free}, and it is a multi-graph if there exist multiple edges connecting the same vertex. The edges $ \mathsf{e} $ and $ \mathsf{e}' $ are said to be incident if they have a common vertex. A path from $ \mathsf{v} $ to $\mathsf{v}' $ or an $ \mathsf{v}-\mathsf{v}' $ path in a $ \mathsf{G} $ is a (loop-free) finite alternating sequence
$
\mathsf{v}=\mathsf{v}_{0},\mathsf{e}_{1},\mathsf{v}_{1},\mathsf{e}_{2},\mathsf{v}_{3},\mathsf{e}_{3},\ldots,\mathsf{e}_{n-1}, \mathsf{v}_{n-1},\mathsf{e}_{n},\mathsf{e}_{n}=\mathsf{v}'
$
of vertices and edges from $ \mathsf{G} $, starting at vertex $ \mathsf{v} $ and ending at vertex $ \mathsf{v}'$ and involving the $ n $ edges $ \mathsf{e}_{k}=(\mathsf{e}_{k-1},\mathsf{e}_{k}) $, where $ 1\leq k\leq n $. The length of this path is $ n $, the number of edges in the path. A path is said to be closed, or is called a cycle, if the start and end of the path coincide. A directed graph $ \mathsf{G} $ is called strongly connected if there is a path between any two distinct vertices of $ \mathsf{G} $. The line graph $L(\mathsf{G})$ of a directed graph $\mathsf{G}$ is the graph obtained
from $\mathsf{G}$ by exchanging the role of the vertices and edges.

Here and in the following, we consider an infinite connected metric graph $\mathsf{G}=(\mathsf{V},\mathsf{E})$ and a flow on it (the latter is described by the boundary controlled vertex delay problem $\mathsf{(vdp)}$). Each edge is normalized so as to be identified with the interval $[0,1]$ and parameterized them contrary to the direction of the flow of material on them, i.e., the material flows from $1$ to $0$. The topology of the graph $\mathsf{G}$ is described by the incidence matrix $\mathcal{I}=\mathcal{I}^{+}-\mathcal{I}^{-}$, where $\mathcal{I}^{-}$ and $\mathcal{I}^{+}$ are \emph{the outgoing incidence} and the \emph{incoming incidence matrices} of $\mathsf{G}$ having entries
\begin{align*}
\textsf{i}^{-}_{ij}:=
\begin{cases}
1,\quad \text{if  } \mathsf{v}_{i}=e_{j}(1),
\\
0,\quad \text{if not,}
\end{cases}
\quad\textsf{i}^{+}_{ij}:=
\begin{cases}
1, \quad \text{if  } \mathsf{v}_{i}=e_{j}(0),
\\
0, \quad \text{if not,}
\end{cases}
\end{align*}
respectively. In particular, we define the outset (resp. inset) of a vertex $ \mathsf{v}_{i} $ as
\begin{align*}
& \text{out}(\mathsf{v}_{i}):=\left\{e_{j} \in \mathsf{E}\,\vert \, \textsf{i}^{-}_{ij}=1 \right\},\qquad\left(\text{resp}.\; \text{in}(\mathsf{v}_{i}):=\left\{e_{j} \in \mathsf{E}\,\vert\,  \textsf{i}^{+}_{ij}=1 \right\}\right).
\end{align*}
The infinite graph $\mathsf{G}$ is called \emph{outgoing or incoming locally finite} if for all $ \mathsf{v}\in \mathsf{V} $ there is $ \delta_{\mathsf{v}}>0 $ such that
$$
\#\text{out}(\mathsf{v})\leq \delta_{\mathsf{v}} \; \; \text{or}\; \; \#\text{in}(\mathsf{v}) \leq \delta_{\mathsf{v}},
$$
and locally finite if it is both outgoing or incoming locally finite \cite{Delio1}. Replacing $1$ by $\textsf{w}_{ij}\geq 0$ in the definition of $\textsf{i}^{-}_{ij}$, we obtain the so-called weighted outgoing incidence matrix
\begin{align*}
\mathcal{I}^{-}_{\mathsf{w}}:=(\mathsf{w}_{ij})_{\mathbb{N}\times \mathbb{N}}.
\end{align*}
In this cases, $\mathsf{G}$ is called a weighted graph and its topology is described via weighted adjacency matrices, for example, the matrix $\mathbb{B}:=(\mathcal{I}_{w}^{-})^{\top}\mathcal{I}^{+}$ is called \emph{the (transposed) weighted adjacency matrix} for the line graph $L(\mathsf{G})$.
With this, one can see that the boundary condition (BC) in the undelayed case, i.e., $ L_{k}=0 $ with $u\equiv 0$ can be written as
\begin{align}\label{S1.1}
\mathcal{I}^{-}c(1)\left(\begin{smallmatrix}
z_{1}(t,1)\\
z_{2}(t,1)\\
\vdots
\end{smallmatrix}\right)=\mathbb{B}c(0)\left(\begin{smallmatrix}
z_{1}(t,0)\\
z_{2}(t,0)\\
\vdots\\
\end{smallmatrix}\right),
\end{align}
where $c(1):=\text{diag }(c_j(1))_{j\in \mathbb{N}}$ and $c(0):=\text{diag }(c_j(0))_{j\in \mathbb{N}}$. If, in addition, we assume that the weights $\textsf{w}_{ij}$ satisfies
\begin{align}\label{S1.2}
\sum_{j\in \mathbb{N}}\textsf{w}_{ij}=1,\; \forall\, i\in \mathbb{N},
\end{align}
then (BC) exhibits standard \emph{Kirchhoff conditions} and the matrices $\mathbb{B}$ are column stochastic.

\section{A concise background on feedback theory of infinite dimensional linear systems}\label{S.2} In this section $\mathscr{X},\mathscr{Z}$ and $\mathscr{U}$ are Banach spaces such that $\mathscr{Z}\subset \mathscr{X}$ with continuous and dense embedding. Throughout this section $\mathscr{A}_m:\mathscr{Z}\to \mathscr{X}$ is a closed linear operator and $\mathscr{G}:\mathscr{Z}\to \mathscr{U}$ is linear surjective operator. We assume that the operator
\begin{align*}
\mathscr{A}:=(\mathscr{A}_m)_{|D(\mathscr{A})}\quad\text{with}\quad D(\mathscr{A}):=\ker\mathscr{G}
\end{align*}
generates a strongly continuous semigroup $\mathscr{T}:=(\mathscr{T}(t))_{t\ge 0}$ on $\mathscr{X}$.

Consider the observed linear system
\begin{align}\label{observed-system}
\begin{cases}
\dot{w}(t)= \mathscr{A}_mw(t),\quad w(0)=w^0, & t>0,\cr \mathscr{G}w(t)=0,& t\ge 0,\cr y(t)=\mathscr{M}w(t), & t\ge 0,
\end{cases}
\end{align}
where $\mathscr{M}:\mathscr{Z}\to\mathscr{U}$ is a linear operator not necessarily closed or closeable. For initial conditions $w^0\in D(\mathscr{A}),$ we have $w(t)=\mathscr{T}(t)w^0$ and then the observation function $y(t)=\mathscr{M}\mathscr{T}(t)w^0$ is well defined for any $t\ge 0$. This is because the domain $D(\mathscr{A})$ is stable by the semigroup $\mathscr{T}$. However, for the instance it is not clear how to define $y(t;w^0)$ for arbitray $w^0\in\mathscr{X}$. We say that the system \eqref{observed-system} is well-posed if the output function $t\mapsto y(t)$ can be extended to a function (denoted by the same symbol) $y\in L^p_{loc}([0,+\infty),\mathscr{U})$ such that
\begin{align*}
\|y(\cdot;w^0)\|_{L^p([0,\al],\mathscr{U})}\le \ga \|w^0\|_{\mathscr{X}}\qquad (w^0\in\mathscr{X})
\end{align*}
for any $\al>0$ and some constant $\ga:=\ga(\al)>0$. Next we will introduce a condition on the operator $\mathscr{C}$ guaranteing the well-posedness of the system \eqref{observed-system}. To this end, we define
\begin{align*}
C:=\mathscr{M}_{|D(C)}\quad\text{with}\quad D(C):=D(\mathscr{A}).
\end{align*}
\begin{definition}\label{2.10}
	The operator $ C$ is called an admissible observation operator for $\mathscr{A} $ if for some (hence all) $\al>0$, there exists a constant $\ga:=\ga(\al)>0$ such that
	\begin{align}\label{admi}
	\int_{0}^{\al} \Vert C \mathscr{T}(t)x\Vert^{p}dt\leqslant \gamma^{p}\Vert x \Vert^{p}\qquad (x\in D(\mathscr{A})).
	\end{align}
	We also say that $(C,\mathscr{A})$ is admissible.
\end{definition}
If $(C,\mathscr{A})$ is admissible, then the following
\begin{align*}
\Psi: D(\mathscr{A})\to L^p_{loc}([0,+\infty),\mathscr{U}),\quad x\mapsto \Psi x= C\mathscr{T}(\cdot)x
\end{align*}
is well defined and has a bounded extension to each $\tilde{\Psi}:\mathscr{X}\to L^p([0,\al],\mathscr{U})$ for any $\al>0$. In this case, the system \eqref{observed-system} is well-posed and the extension of its output function is given by $y=\tilde{\Psi}w^0$ for any $w^0\in\mathscr{X}$. We refer to \cite[chap.3]{TW} and \cite{WO} for more details on admissible observation operators.

In order to given a representation to the extended output function, Weiss \cite{WO} introduced the following Yosida extension of $C$ for $\mathscr{A},$
\begin{align}\label{cgama}
\begin{split}
D(C_{\Lambda})&:=\{x\in \mathscr{X}:\lim_{\lambda\to+\infty} C\lambda R(\lambda,\mathscr{A})x\;\text{exists}\}\cr
C_{\Lambda}x& :=\lim_{\lambda \to+ \infty} C\la  R(\la ,\mathscr{A})x,\qquad x\in D(C_{\Lambda}).
\end{split}
\end{align}
It is shown in \cite{WO} that if $(C,\mathscr{A})$  then for any $x\in \mathscr{X},$ $\mathscr{T}(t)x\in D(C_{\Lambda})$ for a.e. $t>0,$ and the extended output function of the system \eqref{observed-system} is represented as  $y(t;w^0)=C_\Lambda\mathscr{T}(t)w^0$ for initial conditions $w^0\in\mathscr{X}$ and a.e. $t>0$.

Consider the controlled equation
\begin{align}\label{controlled-system}
\begin{cases}
\dot{w}(t)= \mathscr{A}_mw(t),\quad w(0)=w^0, & t>0,\cr \mathscr{G}w(t)=u(t),& t\ge 0,
\end{cases}
\end{align}
where $u:[0,+\infty)\to \mathscr{U}$ is the boundary control function. According to \cite{Gr}, for any $\mu\in\rho(\mathscr{A}),$ the following inverse
\begin{align*}
\mathscr{D}_\mu:=\left(\mathscr{G}_{|\ker(\mu-\mathscr{A}_m)}\right)^{-1}\in\calL(\mathscr{U},\mathscr{X})
\end{align*}
exists and called the Dirichlet operator associated with the boundary control system \eqref{controlled-system}. Now we define the control operator
\begin{align}\label{Bc}
\mathscr{B}=(\mu -\mathscr{A}_{-1})\mathscr{D}_{\mu}\in\calL(\mathscr{U},\mathscr{X}_{-1}), \qquad \mu \in\rho(\mathscr{A}),
\end{align}
where $\mathscr{X}_{-1}$ is the extrapolation space associated with $\mathscr{A}$ and $\mathscr{X},$ and $\mathscr{A}_{-1}:\mathscr{X}\to \mathscr{X}_{-1}$ is the generator of the extrapolation semigroup $\mathscr{T}_{-1}:=(\mathscr{T}_{-1}(t))_{t\ge 0}$ extension of the semigroup $\mathscr{T}$ to $\mathscr{X}_{-1}$. We note that $\mathscr{B}$ does not depends on $\mu,$ due to the resolvent equation. The system \eqref{controlled-system} is reformulated as the following distributed one
\begin{align}\label{2.3}
\dot{w}(t) =\mathscr{A}_{-1}w(t)+\mathscr{B}u(t),\quad w(0) = w^{0},\quad t>0.
\end{align}
If we set
\begin{align*}
\Phi_{t}u:=\int_{0}^{t}\mathscr{T}_{-1}(t -s)\mathscr{B}u(s)ds,
\end{align*}
then the integral solution of the equation \eqref{2.3} is given by $w(t)=\mathscr{T}(t)w^0+\Phi_t u$ for any $t\ge 0,\;w^0\in\mathscr{X}$ and $u\in L^p([0,+\infty),\mathscr{U})$. Observe that $w(t)\in \mathscr{X}_{-1}$ for any $t\ge 0$. The following definition is due to Weiss \cite{WC}.
\begin{definition}\label{2.6}
	We say that the operator  $ \mathscr{B}$ is an admissible control operator for $\mathscr{A}$, if for some $\tau >0 $, we have $ \Phi_{\tau}u\in \mathscr{X}$ for any $u\in L^p([0,+\infty),\mathscr{U})$. In this case, we also say that $(\mathscr{A},\mathscr{B})$ is admissible.
\end{definition}
According to \cite{WC} if $(\mathscr{A},\mathscr{B})$ is admissible, then for any $t\ge 0,$
\begin{align*}
\Phi_t \in \calL(L^p([0,+\infty),\mathscr{U}),\mathscr{X}).
\end{align*}
In this case the solution of \eqref{2.6} satisfies $w(t)\in \mathscr{X}$
for any $t\ge 0,$ $w^0\in\mathscr{X}$ and $u\in L^p([0,+\infty),\mathscr{U})$.

We now consider the input-output system
\begin{align}\label{input-output-system}
\begin{cases}
\dot{w}(t)= \mathscr{A}_mw(t),\quad w(0)=w^0, & t>0,\cr \mathscr{G}w(t)=u(t),& t\ge 0,\cr y(t)=\mathscr{M}w(t),& t\ge 0.
\end{cases}
\end{align}
In the same spirit as for the above two  reformulations, we transform the system \eqref{input-output-system} to the following distributed input-output linear system
\begin{align}\label{DS-input-output}
\begin{cases}
\dot{w}(t) =\mathscr{A}_{-1}w(t)+\mathscr{B}u(t),\quad w(0) = w^{0},& t>0,\cr y(t)=Cw(t),& t\ge 0.\end{cases}
\end{align}
The following definition can be found in \cite{Sa}, \cite{St} and \cite{WR}.
\begin{definition}\label{well-posed-ABC}
	We say that the system \eqref{DS-input-output} (or the triple $(\mathscr{A},\mathscr{B},C)$) is well-posed if and only if $(\mathscr{A},\mathscr{B})$ is admissible and the output function $y$ of the system is extended to a function $y\in L^p_{loc}([0,+\infty),\mathscr{U})$ such that
	\begin{align}\label{y-estima}
	\|y(\cdot;w^0,u)\|_{L^p([0,\al],\mathscr{U})} \le \kappa \left(\|w^0\|_{\mathscr{X}}+\|u\|_{L^p([0,\al],\mathscr{U})}\right)
	\end{align}
	for any $w^0\in\mathscr{X}$, any $\al>0$ and any $u\in L^p([0,\al],\mathscr{U}),$ where $\kappa:=\kappa(\al)>0$ is a constant.
\end{definition}
An necessarily condition for the well-posedness of the triple $(\mathscr{A},\mathscr{B},C)$ is that $(C,\mathscr{A})$ is admissible. In fact, it suffice to  choose $w^0\in D(\mathscr{A})$ and $u=0$ and use the estimate \eqref{y-estima}. In order to shed more light on the  well-posedness of the triple $(\mathscr{A},\mathscr{B},C)$, we define the following dense spaces
\begin{align*}
W^{1,p}_{0,\al}(\mathscr{U}):=\left\{u\in W^{1,p}([0,\al],\mathscr{U}):u(0)=0\right\},\qquad \al>0.
\end{align*}
Now assume that $(\mathscr{A},\mathscr{B})$ is admissible and (without loss of generality) $0\in\rho(\mathscr{A})$. An integration by parts shows that
\begin{align*}
\Phi_{t}u=\mathscr{D}_0u(t)-R(0,\mathscr{A})\Phi_t \dot{u}\in\mathscr{Z}
\end{align*}
for any $t\in  [0,\al],$ and  $u\in W^{1,p}_{0,\al}(\mathscr{U})$. This allows us to introduce the following map
\begin{align*}
(\mathbb{F}u)(t):=\mathscr{M}\Phi_{t}u,\qquad u\in W^{1,p}_{0}([0,\al],\mathscr{U}),\;a.e.\;t\in [0,\al].
\end{align*}
For any $\al>0$ and $(w^0,u)\in D(\mathscr{A})\times W^{1,p}_{0}([0,\al],\mathscr{U}),$ the output function of the system \eqref{DS-input-output} is given by
\begin{align*}
y(t;w^0,u)=C\mathscr{T}(t)w^0+(\F u)(t),\qquad t\in [0,\al].
\end{align*}
By using this expression together with Definition \ref{well-posed-ABC}, one easily prove the following result.
\begin{proposition} \label{caracterisation-wellposedABC}
	The triple $(\mathscr{A},\mathscr{B},C)$ is well-posed if and only if the following assertions hold:
	\begin{enumerate}
		\item $(\mathscr{A},\mathscr{B})$ and $(C,\mathscr{A})$ are admissible.
		\item for $\al>0,$ there exists a constant $\kappa:=\kappa(\al)>0$ such that
		\begin{align}\label{est-F}
		\|\F u\|_{L^p([0,\al],\mathscr{U})}\le \kappa \| u\|_{L^p([0,\al],\mathscr{U})}
		\end{align}
		for any $u\in  W^{1,p}_{0}([0,\al],\mathscr{U})$.
	\end{enumerate}
	In this case, the operator $\F$ is extended to $\F\in \calL(L^p([0,\al],\mathscr{U}))$ and the extended output function $y$ of the system \eqref{DS-input-output} satisfies
	\begin{align*}
	y(t;w^0,u)=C_\Lambda \mathscr{T}(t)w^0+(\F u)(t)
	\end{align*}
	for any $w^0\in\mathscr{X},$ $u\in L^p_{loc}([0,+\infty),\mathscr{U})$ and a.e. $t>0$. The operator $\F$ is called the extended input-output control operator.
\end{proposition}
In order to give a complete representation of the extended output function $y$, we need the following subclass of well-posed linear systems introduced in \cite{WR}.
\begin{definition}\label{Reg}
	A well-posed triple $(\mathscr{A},\mathscr{B},C)$ with extended input-output control operator $\F$ is called regular (with feedthrough zero) if for any $ v\in \mathscr{U} $, the following limit exists in $\mathscr{U}$.
	\begin{align*}
	\lim_{\tau\longmapsto 0}\frac{1}{\tau}\int_{0}^{\tau}(\mathbb{F}(\1_{\mathbb{R}_{+}} \cdot v))(\sigma)d\sigma=0.
	\end{align*}
\end{definition}
The following result is proved in \cite{Sa} and \cite{WR}.
\begin{theorem}\label{weiss-rep}
	Let $(\mathscr{A},\mathscr{B},C)$ be a regular triple with extended input-output control operator $\F$. Then $\Phi_t u\in D(C_\Lambda)$ and $(\F u)(t)=C_\Lambda \Phi_t u$ for any $u\in L^p([0,+\infty),\mathscr{U})$ and a.e. $t\ge 0$. In particular, we have the state trajectory and the output function of the system \eqref{DS-input-output} satisfy $w(t;w^0,u)\in D(C_\Lambda)$ and $y(t;w^0,u)=C_\Lambda w(t;w^0,u)$ for any $w^0\in\mathscr{X},$ $u\in L^p_{loc}([0,+\infty),\mathscr{U})$ and a.e. $t> 0$.
\end{theorem}
In what follow we recall the concept of feedback theory for regular systems.
\begin{definition}\label{admi-feed}
	Let $(\mathscr{A},\mathscr{B},C)$ be a well-posed triple with extended input-output control operator $\F$. We say that the identity operator $I:\mathscr{U}\to\mathscr{U}$ is an admissible feedback operator for $(\mathscr{A},\mathscr{B},C)$ if for some (hence all) $\al>0$, the operator $I-\F$ has a uniformly bounded inverse in $\calL(L^p([0,\al],\mathscr{U}))$.
\end{definition}
We now consider the following operator
\begin{align}\label{2.15}
\mathfrak{A}=\mathscr{A}_{m},  \qquad D(\mathfrak{A})=\lbrace x\in \mathscr{Z}:\mathscr{G}x=\mathscr{M}x\rbrace,
\end{align}
where the operators  $\mathscr{A}_{m},$ $\mathscr{G}$ and $\mathscr{M}$ are as above.  The following result is proved in \cite[Theorem 4.1]{HMR}
\begin{theorem}\label{T.1}
	Assume that $(\mathscr{A},\mathscr{B},C)$ is a regular triple with the identity operator $ I:\mathscr{U}\to \mathscr{U} $ as an admissible feedback. Then the operator $\mathfrak{A}$ coincides with the following one
	\begin{align*}
	\mathscr{A}^{cl}&=\mathscr{A}_{-1}+\mathscr{B}C_\Lambda\cr D(\mathscr{A}^{cl})&=\left\{x\in D(C_\Lambda):(\mathscr{A}_{-1}+\mathscr{B}C_\Lambda)x\in \mathscr{X}\right\},
	\end{align*}
	which generates a strongly continuous semigroup $\mathfrak{T}:=(\mathfrak{T}(t))_{t\geq 0}$ on $\mathscr{X}$. Furthermore, we have $\mathfrak{T}(s)x\in D(C_\Lambda)$ for any $x\in \mathscr{X}$ and a.e. $s>0$, and
	\begin{align*}
	\mathfrak{T}(t)x=\mathscr{T}(t)x+\int^t_0 \mathscr{T}_{-1}(t-s)\mathscr{B}C_\Lambda \mathfrak{T}(s)x\,ds
	\end{align*}
	for any $t\ge 0$ and $x\in \mathscr{X}$. In addition, for $\mu\in\rho(\mathscr{A})$, we have
	\begin{align*}
	\mu\in\rho(\mathfrak{A}) \Leftrightarrow 1\in \rho(\mathscr{D}_{\mu}\mathscr{M}) \Leftrightarrow 1\in \rho(\mathscr{M}\mathscr{D}_{\mu}).
	\end{align*}
	In this case,
	\begin{align}\label{lap}
	\begin{split}
	R(\mu,\mathfrak{A})&=(I-\mathscr{D}_{\mu}\mathscr{M})^{-1}R(\mu,\mathscr{A})\\\nonumber
	&=\left(I+\mathscr{D}_{\mu}(I-\mathscr{M}\mathscr{D}_{\mu})^{-1}\mathscr{M}\right)R(\mu,\mathscr{A}).
	\end{split}
	\end{align}
\end{theorem}
We end this section by recalling the well-posedness of the following  inhomogeneous boundary problem
\begin{equation}\label{Sy.2}
\begin{cases}
\dot{w}(t)= \mathscr{A}_{m} w(t),\quad w(0)=w^0,&  t> 0,\cr
\mathscr{G} w(t)= \mathscr{M}w(t)+f(t),& t\geq 0,
\end{cases}
\end{equation}
where $f:\R_+\to\mathscr{U}$ is a locally $p$-integrable function. The following result is proved in \cite[Theorem 4.3]{HMR}.
\begin{theorem}\label{nonhomo-bound}
	Assume that $(\mathscr{A},\mathscr{B},C)$ is a regular triple with the identity operator $ I:\mathscr{U}\to \mathscr{U} $ as an admissible feedback. Then for any initial conditions $w^0\in \mathscr{X}$ and any $f\in L^p_{loc}([0,+\infty),\mathscr{U}),$ the system \eqref{Sy.2} has a unique mild solution $w:[0,+\infty)\to\mathscr{X}$ given by
	\begin{align*}
	w(t)=\mathfrak{T}(t)w^0+\int^t_0 \mathfrak{T}_{-1}(t-s)\mathscr{B}f(s)\,ds,\qquad t\ge 0,
	\end{align*}
	where $(\mathfrak{T}(t))_{t\ge 0}$ is the strongly continuous semigroup generated by the operator $\mathfrak{A}$ defined in \eqref{2.15}.
\end{theorem}

\section{Well-posedness of the vertex delay system } \label{S.3}
In this section we are concerned with the existence and uniqueness of the solutions of the vertex delay system $\mathsf{(vdp)}$. Inspirited from \cite{HI}, \cite[Section 5]{HMR}, we will use product spaces to reformulate the system $\mathsf{(vdp)}$ as a free-delay boundary linear system. This allows us to  use  Theorem \ref{T.1} and Theorem \ref{nonhomo-bound} to prove the well-posedness of $\mathsf{(vdp)}$.
\subsection{The boundary delay network equations as abstract delay boundary value problem}\label{sec3-sub1}
The object of this section is to introduce notation and appropriate operators in order to rewrite the network delay equation $\mathsf{(vdp)}$ in abstract way. Before that, we first impose the following conditions on the coefficients appearing in the system $\mathsf{(vdp)}$.
\begin{mainassumptions}
	Assume that:
	\begin{itemize}
		\item[{\bf (A1)}]  $ c_j(\cdot),q_j(\cdot)\in L^{\infty}([0,1]) $ such that $c_j(x)\ge \ga_1 , q_j(x)\le \ga_2$ for $j\in \mathbb{N},\,x\in [0,1]$ and some constants $\ga_1 >0,\,\ga_2\in \mathbb{R}$.
		\item[{\bf (A2)}] Each vertex has only finitely many outgoing edges.
		\item[{\bf (A3)}] The weights $\textsf{w}_{ij}$ satisfies \eqref{S1.2}.
		\item[{\bf (A4)}] For any $i\in\N$, $\eta_i:[-r,0]\to\R$ is a function of bounded variations such that $|\eta_i|([-\varepsilon,0])\to 0$ as $\varepsilon\to 0$, and for any $\theta\in [-r,0],$ $\eta(\theta)={\rm diag}(\eta_i(\theta))_{i\in\N}\in\ell^1$.
	\end{itemize}
\end{mainassumptions}
\begin{remark}
	In the above, {\bf (A1)} specifies the transport process along each edges, {\bf (A2)} is equivalent the statement that $\mathsf{G}$ is an infinite outgoing locally finite graph, while {\bf (A3)} enforces the conservation of flows at vertices: the amount of material flowing into a vertex $\mathsf{v}_i$ must equal the amount that flows out from this vertex, this is so far for all vertices except the source.
\end{remark}
In view of {\bf (A1)}, we define
\begin{align}\label{S3.functions}
\tau_{j}(x_1,x_2):=\int_{x_1}^{x_2}\dfrac{dx}{c_{j}(x)},\qquad
\xi_{j}(x_1,x_2):=\int_{x_1}^{x_2} \dfrac{q_{j}(x)}{c_{j}(x)}dx,
\end{align}
for every $j\in \mathbb{N}$ and $x_1,x_2\in [0,1]$. These functions define the time required to move between two locations on the edge $\mathsf{e}_j $ with the velocity $c_j(\cdot)$ and the rate of mass gained or lost during this movement resulting from the factor $ q_j(\cdot) $, respectively. Moreover, {\bf (A2)} and {\bf (A3)} yield that the outgoing incident matrix $\mathcal{I}^{-}$ is bounded on $\ell^1$ and the weighted (transposed) adjacency matrix $\mathbb{B}$ satisfy $\Vert\mathbb{B}\Vert_{\ell^1}=1$, respectively.

In the rest of this paper,  $ p>1 $ and $r>0$ are real numbers. Define Banach spaces by $(X,\|\cdot\|_X)$ and $(Y,\|\cdot\|_Y)$ by
\begin{align*}
& X:=L^{p}([0,1],\ell^{1}),\quad \Vert f \Vert^{p}_X :=\int_{0}^{1}\Vert f(x)\Vert^{p}_{\ell^{1}} dx,\cr
& Y:=L^{p}([-r,0],X), \quad \|g \|^{p}_Y :=\int_{-r}^{0}\int^1_0 \|g(\theta,x)\|^{p} dx d\theta.
\end{align*}
In order to rewrite in abstract way, the following differential operators are needed
\begin{align}\label{S1.A_m}
A_{m}g:= c(\cdot)\partial_x g+q(\cdot)g,\quad
D(A_{m}):= W^{1,p}([0,1];\ell^{1})
\end{align}
and
\begin{align}\label{S1.Qm}
Q_{m}\varphi =\partial_{\theta}\varphi, \;\; D(Q_{m})= W^{1,p}([-r,0],X),
\end{align}
where $c={\rm diag}(c_j)_{j\in\N}$ and  $q={\rm diag}(q_j)_{j\in\N}$.
Moreover, for $\varphi=(\varphi_k)_{k\in\N}\in D(Q_m),$ we select $L(\varphi)=(L_k(\varphi_k))_{k\in\N},$ where
\begin{align}\label{S1.3}
L_k(\varphi_k)=\int_{0}^{1}\int_{-r}^{0}d[\eta_k(\theta)]c_k(x)\varphi_k(\theta,x)dx,\;k\in\N.
\end{align}
The following operators are required in our transformation
\begin{align}\label{3.13}
\begin{split}
&\G g:=	\mathcal{I}^{-} c(1)g(1),\quad  g\in W^{1,p}([0,1];\ell^{1}),\\
&\M g:=\mathbb{B}c(0)g(0),\quad  g\in W^{1,p}([0,1];\ell^{1}),\\
& \L\varphi:=\mathbb{B}L\varphi, \quad \varphi\in D(Q_m),\cr
&K:=(\mathsf{b}_{il})_{i,l\in\mathbb{N}\times \{1,\ldots,N\}}.
\end{split}
\end{align}
We set $z:=(z_j)_{j\in\N}:[-r,+\infty)\to X$, where $z_j$ are given in the system $\mathsf{(vdp)}$. The history function of $z$ at the time $t\geq 0$ is the function $z_t:[-r,0]\to X$ defined by $z_t(\theta,x)=z(t+\theta,x)$ for any $\theta\in [-r,0]$. The function $z_0=\varphi:=(\varphi_{j})_{j\in \mathbb{N}}$ is called the initial history function of $z(\cdot,\cdot)$. Using the previous notation, we can rewrite the system $\mathsf{(vdp)}$ as the following abstract delay boundary value problem on $X,$
\begin{align}\label{GML}
\begin{cases}
\dot{z}(t)=A_m z(t),& t\ge 0,\cr \G z(t)=\M z(t)+\L z_t+Ku(t),& t\ge 0,\cr z(0)=g,\quad z(t)=\varphi(t),& t\in [-r,0],
\end{cases}
\end{align}
where the control function $u:[0,+\infty)\to U:=\C^N$ is given by $u=(u_1,\cdots, u_N)$.

As we work with delay system then is more convenient to introduce the new state space
\begin{align*}
\mathscr{X} :=X\times Y\quad\text{with norm}\quad \left\|(\begin{smallmatrix} g\\ \varphi\end{smallmatrix})\right\|:=\|g\|_X+\|\varphi\|_Y.
\end{align*}
Furthermore, we introduce the spaces  \begin{align*}\mathscr{Z}:= D(A_{m})\times D(Q_{m}),\quad \mathscr{U}:=\ell^1\times X,\end{align*} and the following matrix operator
\begin{align}\label{S1.GM}
\begin{split}
& \mathscr{A}_{m}:=\begin{pmatrix}A_m& 0\\ 0& Q_m\end{pmatrix}: \mathscr{Z}\to\mathscr{X},\qquad
\mathscr{G}:=\begin{pmatrix}
\G & 0 \\
0 & \delta_{0} \\
\end{pmatrix},\qquad \mathscr{Z} \to \mathscr{U}\\
&\mathscr{M}:=\begin{pmatrix}
\M &\L \\
I_{X} & 0 \\
\end{pmatrix},\qquad \mathscr{Z} \to \mathscr{U}, \qquad \mathfrak{B}:=(\begin{smallmatrix}K\\0\end{smallmatrix}):U\to \mathscr{U}.
\end{split}
\end{align}
Now by selecting the new state
\begin{align*}
w(t)=\left(\begin{smallmatrix} z(t)\\ z_{t}\end{smallmatrix}\right),\qquad t\ge 0,
\end{align*}
the delay boundary control problem \eqref{GML} (hence the delay network system $\mathsf{(vdp)}$) is reformulated in $\mathscr{X}$ as the following free delay perturbed boundary control system
\begin{align}\label{Sy.6}
\left\lbrace
\begin{array}{lll}
\dot{w}(t)&=&\mathscr{A}_{m}w(t),\qquad\qquad\; t> 0, \\
w(0)&=&w^0,\\
\mathscr{G} w(t)&=&\mathscr{M} w(t)+\mathfrak{B}u(t),\quad t> 0,
\end{array}
\right.
\end{align}

\subsection{Solving a free delay boundary network equation}\label{sec3-sub2}
In this section we will study the well-pposedness of the system \eqref{Sy.6} (hence the netwok delay equation $\mathsf{(vdp)}$). In fact, as in Section \ref{S.2}, to the system \eqref{Sy.6}, we associate the following operator
\begin{align}\label{frakA}
\mathfrak{A}:=\mathscr{A}_{m},\qquad 
D(\mathfrak{A}):=  \big\lbrace \left(\begin{smallmatrix}g\\ \varphi\end{smallmatrix}\right) \in D(\mathscr{A}_{m}): \mathscr{G}(\begin{smallmatrix}g\\ \varphi\end{smallmatrix})=\mathscr{M}(\begin{smallmatrix}g\\ \varphi\end{smallmatrix})  \big\rbrace.
\end{align}
where $\mathscr{A_m},\mathscr{G}$ and $\mathscr{M}$ are defined in \eqref{S1.GM}.

\begin{remark}\label{important-rem}
	We mention that
	\begin{align*}
	D(\mathfrak{A})& =
	\big\lbrace \left(\begin{smallmatrix}g\\ \varphi\end{smallmatrix}\right) \in D(\mathscr{A}_{m}) : \qquad Gg=Mg+ \mathrm{L}_{\mathbb{B}}\varphi,\; g=\varphi(0)\big\rbrace\nonumber,
	\end{align*} By the same argument as in \cite[Proposition 3.1]{Dorn-infinie}, for $\left(\begin{smallmatrix}g\\ \varphi\end{smallmatrix}\right) \in \mathscr{Z}$,
	\begin{align*}
	\left(\begin{smallmatrix}g\\ \varphi\end{smallmatrix}\right)\in D(\mathfrak{A})\Leftrightarrow g(1)=M g+ R \varphi,
	\end{align*}
	where
	\begin{align}\label{walid}
	M:=c^{-1}(1) \M= c^{-1}(1)\B c(0)\delta_0,\qquad R :=c^{-1}(1)\L=c^{-1}(1)\B L.
	\end{align}
	Thus
	\begin{align*}
	D(\mathfrak{A}):=  \big\lbrace \left(\begin{smallmatrix}g\\ \varphi\end{smallmatrix}\right) \in D(\mathscr{A}_{m}): \widetilde{\mathscr{G}}(\begin{smallmatrix}g\\ \varphi\end{smallmatrix})=\widetilde{\mathscr{M}}(\begin{smallmatrix}g\\ \varphi\end{smallmatrix})  \big\rbrace\nonumber,
	\end{align*}
	where
	\begin{align}\label{S1.GM-tilde}
	& \widetilde{\mathscr{G}}:=\begin{pmatrix}
	\delta_1 & 0 \\
	0 & \delta_{0} \\
	\end{pmatrix},\quad \mathscr{Z} \to \mathscr{U}\qquad
	&\widetilde{\mathscr{M}}:=\begin{pmatrix}
	M &R \\
	I_{X} & 0 \\
	\end{pmatrix},\quad \mathscr{Z} \to \mathscr{U}.
	\end{align}
\end{remark}

Next we will prove that $\mathfrak{A}$ is a generator on $\mathscr{X}$. To this end we need some preliminaries.  Define the operators \begin{align*}
&G:D(A_m)\to\ell^1,\quad g\mapsto Gg=g(1), \cr & A=(A_m)_{|D(A)}\quad\text{with}\quad D(A):= \ker G\cr & Q=(Q_m)_{|D(Q)}\quad\text{with}\quad D(Q)=\ker\delta_0.
\end{align*}
A similar argument as in \cite[Lemma 3.3]{TE} shows that under the condition {\bf (A1)} the operator $A$ generates a strongly continuous nilpotent semigroup $(T(t))_{t\geqslant 0} $ on $X$ given by
\begin{align}\label{3.17}
(T(t)g)_{j}(x)=
\begin{cases}
e^{\xi_j(x, \tilde{s}(t))}g_{j}(\tilde{s}_j(t)), & {\rm if } \; t \leq \tau_{j}(x,1),
\\
0,& {\rm if\; not }.
\end{cases}
\end{align}
for all $x\in [0,1]$ and $j\in \mathbb{N}$. Here $\tilde{s}_j(t) $ is a continuous function defined by $\tau_{j}(x,\tilde{s}_j(t) )=t$, which determines the location where the flow moves on the edge $ e_{j} $ from the point $ x$ during time $ t\leq \tau_{j}(x,1) $. In addition, the resolvent set of $A$ is $\rho(A)=\C$. On the other hand, it is well known that (see e.g. \cite[chap.II]{EN}) the operator $Q$ generates the left shift semigroup $(S(t))_{t\ge 0}$ on $L^p([-r,0],X)$ given by
\begin{align}\label{S-sg}
(S(t)\varphi)(\theta)=\begin{cases} 0,& -t\le \theta\le 0,\cr \varphi(t+\theta),& -r\le \theta\le -t,\end{cases}
\end{align}
for any $t\ge 0$ and $\varphi\in L^p([-r,0],X)$. Moreover, we have $\rho(Q)=\C$.

The proof of the following result can easily obtained by the same computations as in \cite[Theorem 3]{HIR} and \cite[Example 5.2]{HMR}.
\begin{lemma}\label{regular-shift-triple}
	Let the assumptions {\bf  (A1)} to {\bf  (A4)} be satisfied and assume that $\ga_3:=\sup_{i\in\N}\|c_i\|_\infty<\infty$. If we denote by $d_\la$ for $\la\in\C$ the Dirichlet maps associated with $Q_m$ and $\delta_0,$ then $d_\la: X\to \to L^p([-r,0],X))$ is given by \begin{align}\label{Dirichlet-d}(d_\la g)(\theta,x)=e^{\la \theta}g(x),\quad\theta\in [-r,0],\;x\in [0,1].\end{align} Moreover, if $R$ be the delay operator defined in \eqref{walid} and \begin{align}\label{beta}\beta:= (-Q_{-1})d_0:X\to L^p([-r,0],X)_{-1},\end{align} then the triple $(Q,\beta,R)$ is a regular system with control
	\begin{align}\label{control-maps-beta}
	(\Phi^{\beta,Q}_t v)(\theta,x)=\begin{cases} v(t+\theta,x),& -t\le \theta \le 0,\cr 0,& -r\le \theta\le -t,\end{cases}
	\end{align}
	for any $t\ge 0$ and $x\in [0,1]$, and $v\in L^p(\R^+,\ell^1)$.
\end{lemma}
The following result introduce an important regular system required to prove the main results of the rest of the work.
\begin{lemma}\label{L.2}
	Let the assumption {\bf (A1)} and {\bf (A3)} be satisfied.  The Dirichlet operator associated with $A_m$ and $G$ is given by
	\begin{align}\label{Dirichlet-D}
	D_\la =\left({\rm diag}(e^{\xi_j(\cdot,1)-\la \tau_j(\cdot,1)})\right)_{j\in\N}.
	\end{align}
	Moreover we define
	\begin{align}\label{B}
	B:=(\la-A_{-1})D_\la: \ell^1\to X_{-1},\quad\la\in\rho(A),
	\end{align}
	then $(A,B)$ is admissible with control maps $\Phi^{A,B}_t:L^p(\R^+,\ell^1)\to X$ such that for any $j\in\N$
	\begin{align}\label{control-maps-B}
	\left(\Phi^{A,B}_t v\right)_j(x)=\begin{cases} e^{\xi_j(x,1)}v_j(t-\tau_j(x,1)),& t\ge \tau_j(x,1),\cr 0,&\text{if not,}\end{cases}
	\end{align}
	for any $t\ge 0,$ $x\in [0,1]$ and $v\in L^p(\R^+,\ell^1)$. Furthermore, if we assume that $\ga_3:=\sup\{\|c_i\|_\infty:i\in\N\}<\infty$ and $\t_0:=\inf_{ j\in \mathbb{N}}\t_j(0,1)>0$, then triple $(A,B,C)$ is regular, where $C:=M_{|D(A)}$ and $M$ defined in \eqref{walid}.
\end{lemma}
\begin{proof}
	The computation of $D_\la$ can be found e.g. in \cite[p.439]{TE}. On the other hand, by tacking Laplace transform in both side of \eqref{control-maps-B}, on can see that $\widehat{\Phi^{A,B}_\cdot v}(\la)=D_\la \hat{v}(\la)$ for a large $\la>0$, where $\hat{v}$ denotes the Laplace transform of $v$. Now according the injectivity of the Laplace transform and  \eqref{B}, we have
	\begin{align*}
	\int^t_0 T_{-1}(t-s)Bv(s)ds=\Phi^{A,B}_t v,\quad t\ge0,\; v\in L^p(\R^+,\ell^1).
	\end{align*}
	Let us prove that for each $t\ge 0,$ the map $\Phi^{A}_{t}$ is linear bounded from $L^{p}(\R_+,\ell^1)$ to $X$. To this end, we define the following function
	\begin{align}\label{S2.hf}
	h_j(x):=\tau_{j}(x,1),\qquad x\in [0,1],\;j\in\mathbb{N}.
	\end{align}
	Clearly, the function $h(\cdot)$ is continuous on $[0,1]$ and strictly decreasing on $(0,1)$ and hence is invertible. For $t\geq \tau_{j}(x,1)$,  we have $h^{-1}(t)\leq x$. Thus
	\begin{align*}
	& \|\Phi^{A}_{t}v\|_X= \int_{0}^{1}\left(\sum_{j=1}^{+\infty}\left\vert  (\Phi^{A}_{t}v)_{j}(x)\right\vert\right)^p dx\\
	&= \int_{h_j^{-1}(t)}^{1}\left(\sum_{j=1}^{+\infty}\left\vert e^{\xi_{j}(x,1)}v_{i}(t-h_j(x))\right\vert\right)^p dx\\
	&=\int_{0}^{t}\sum_{j=1}^{+\infty}\left\vert e^{\xi_{j}(h_j^{-1}(t-r),1)}v_{i}(r)\right\vert^p  c_j(h_j^{-1}(t-r))dr\cr &\leq \ga_3 e^{p\frac{\ga_2}{\ga_1}}  \int^t_0 \|v(r)\|^p_{\ell^1} dr.
	\end{align*}
	Thus the admissibility of the control operator $B$ for $A$ follows. Let us prove that $(C,A)$ is admissible. In fact, we put $\Upsilon=c^{-1}(1)\B c(0)$ and take  $ g\in D(A) $ and $ \alpha >0 $.  As the matrices $c^{-1}(1)$ and $c(0)$ are diagonal, and $\|\B\|=1,$ then using product matrices properties and the condition {\bf (A1)}, we obtain $\|\Upsilon\|\le \frac{\ga_3}{\ga_1}$. Thus
	\begin{align*}
	\int_{0}^{\alpha}\left\Vert CT(t)g \right\Vert_{\ell^1}^{p}dt
	&=\int_{0}^{\alpha} \left\Vert\Upsilon[T(t)g](0)\right\Vert_{\ell^1}^{p}dt\\
	&\leq \left( \frac{\ga_3}{\ga_1}\right)^p
	\int_{0}^{\min(\tau_{j}(0,1),\alpha)}\left(\sum_{j=1}^{+\infty} \vert e^{\xi_j(0, \tilde{s}(t))}g_{j}(\tilde{s}_j(t)) \vert\right)^{p}dt\\
	&\leq \left( \frac{\ga_3}{\ga_1}\right)^p
	\int_{0}^{1}\left(\sum_{j=1}^{+\infty} \left\vert e^{\xi_j(0,\sigma)}g_{j}(\sigma) \right\vert\right)^{p} \dfrac{1}{c_{j}(\sigma)} d\sigma,
	\end{align*}
	where we used the change of variables $ \sigma= \tilde{s}(t) $ so that
	$$  d\sigma =c_{j}(\tilde{s}(t))dt= c_{j}(\sigma)dt. $$
	The assumptions ({\bf A1}) to ({\bf A3}) imply
	\begin{align*}
	\int_{0}^{\alpha}\left\Vert CT(t)g \right\Vert_{\ell^{1}}^{p}dt \leq
	\frac{\ga_3^p}{\ga_1^{p+1}}  e^{\tfrac{p\ga_2}{\ga_1}}\Vert g\Vert^{p}_X.
	\end{align*}
	It follows that $ C $ is an admissible observation operator for $ A$. As in Section \ref{S.2}, we define the following operators
	\begin{align*}
	(\F^{A} v)(t)= \Upsilon \delta_0 \Phi^{A,B}_t v,\quad v\in W^{1,p}_{0,\al}(\ell^1).
	\end{align*}
	For any $j\in\N$ and  $\al>\tau_j(0,1),$ and $v\in W^{1,p}_{0,\al}(\ell^1),$ we have
	\begin{align*}
	 \|\mathbb{F}^{A}v\|_{L^p([0,\al],\ell^1)}^p &\leq  \|\Upsilon\|^p \int_{\tau_j(0,1)}^{\al}\left\vert\sum_{j=1}^{+\infty} e^{\xi_{j}(0,1)}v_{i}(t-\tau_j(0,1))\right\vert^pdt\cr & \le  \left(\frac{\ga_3}{\ga_1}\right)^p \ga_3 e^{\tfrac{p\ga_2}{\ga_1}} \Vert v\Vert_{L^p([0,\al],\ell^1)}^p,
	\end{align*}
	due to the same calculus as for $\Phi^{A,B}$. Thus by Proposition \ref{caracterisation-wellposedABC}, the triple $(A,B,C)$ is well-posed. In addition, for $\al<\tau_0,$ and $v\in W^{1,p}_{0,\al}(\ell^1),$ we have $(\mathbb{F}^{A}v)(t)=0$ for any $t\in [0,\al]$. In particular, for small $\al>0$, and any $v_0\in \ell^1,$ we have
	\begin{align*}
	\frac{1}{\al}\int^\al_0 (\F^A (\1_{\R^+}\cdot v_0))(\si)d\si=0.
	\end{align*}
	Hence $(A,B,C)$ is regular.
\end{proof}
The first main result of the section is the following:
\begin{theorem}\label{T.2}
	Let the assumptions {\bf  (A1)} to {\bf  (A4)} be satisfied and assume that $\sup_{i\in\N}\|c_i\|_\infty<\infty$ and  $\t_0:=\inf_{ j\in \mathbb{N}}\t_j(0,1)>0$. Then the operator $\mathfrak{A}$ defined by \eqref{frakA} generates a strongly continuous semigroup $\mathfrak{T}:=(\mathfrak{T}(t))_{t\ge 0}$ on $\mathscr{X}$.
\end{theorem}
\begin{proof}
	According to Remark \ref{important-rem}, we have
	\begin{align*}D(\mathfrak{A})=\left\{(\begin{smallmatrix}g\\\varphi\end{smallmatrix})\in\mathscr{Z}: \widetilde{\mathscr{G}}(\begin{smallmatrix}g\\\varphi\end{smallmatrix})=\widetilde{\mathscr{M}}(\begin{smallmatrix}g\\\varphi\end{smallmatrix})\right\}.
	\end{align*}
	Thus we will use Theorem \ref{T.1} to prove that $\mathfrak{A}$ is a generator. In fact,  the operator
	\begin{align*}
	\mathscr{A}=(\mathscr{A}_m)_{|D(\mathscr{A})}\quad\text{with}\quad D(\mathscr{A})=\ker\widetilde{\mathscr{G}}=D(A)\times D(Q)
	\end{align*}
	generates the following diagonal semigroup $\mathscr{T}(t)={\rm diag}(T(t),S(t))$ for  $t\ge 0$, where the semigroups $T(\cdot)$ and $S(\cdot)$ are given by \eqref{3.17} and \eqref{S-sg}, respectively. On the other hand, the Dirichlet maps associated with $\mathscr{A}_m$ and $\mathscr{G}$ are
	$\mathscr{D}_\la= {\rm diag}(D_\la,d_\la)$ for $\la\in\C,$ where $D_\la$ and $d_\la$ are, respectively, given by \eqref{Dirichlet-D} and \eqref{Dirichlet-d}. Now if we define a control operator $\mathscr{B}:=(\la-\mathscr{A}_{-1})\mathscr{D}_\la:\ell^1\times \ell^1\to\mathscr{X}_{-1}$ for $\la\in\C,$ then $\mathscr{B}={\rm diag}(B,\beta),$ where $B$ and $\beta$ are given by \eqref{B} and \eqref{beta}, respectively. Now as $\mathscr{A}$ has a diagonal domain, then $(\mathscr{A},\mathscr{B})$ is well-posed, due to Lemma \ref{L.2} and Lemma \ref{regular-shift-triple}. Moreover, the control maps associated with $(\mathscr{A},\mathscr{B})$ are given by $
	\Phi^{\mathscr{A},\mathscr{B}}_t={\rm diag}(\Phi^{A,B}_t, \Phi^{Q,\beta}_t)$ for any $t\ge 0$. Denote $\mathscr{C}:=\widetilde{\mathscr{M}}$ with domain  $D(\mathscr{C})=D(\mathscr{A})$. We will prove that $(\mathscr{A},\mathscr{B},\mathscr{C})$ is regular. In fact, by combining Lemma \ref{regular-shift-triple} and Lemma \ref{L.2}, one can easily prove that $(\mathscr{C},\mathscr{A})$ is admissible. On the other hand
	\begin{align*}
	\F^{\calA} := \mathscr{C}\Phi^{\mathscr{A},\mathscr{B}}_{\cdot}= \begin{pmatrix}\F^A&\F^Q\\ \Phi^{A,B}_{\cdot}&0\end{pmatrix}
	\end{align*}
	on $W^{1,p}_{0,\al}(\mathscr{U})$. Thus according to the above two lemmas and  Proposition \ref{caracterisation-wellposedABC}, the triple $(\mathscr{A},\mathscr{B},\mathscr{C})$ is well-posed. On the other hand, for small $\al>0$, and constant control $(\begin{smallmatrix}v_1\\v_2\end{smallmatrix})\in\mathscr{U}$, we have
	\begin{align*}
	\left\|\frac{1}{\al}\int^\al_0 (\F^{\mathscr{A}} (\1_{\R^+}\cdot (\begin{smallmatrix}v_1\\v_2\end{smallmatrix})))(\si)d\si\right\|=
	\left\|\frac{1}{\al}\int^\al_0 (\F^Q (\1_{\R^+}\cdot v_2))(\si)d\si\right\|.
	\end{align*}
	Thus $(\mathscr{A},\mathscr{B},\mathscr{C})$ regular due to Lemma \ref{regular-shift-triple}. Let us now show that the identity operator $I_{\mathscr{U}}:{\mathscr{U}}\to{\mathscr{U}}$ is an admissible feedback for the system $(\mathscr{A},\mathscr{B},\mathscr{C})$. For $\al_0<\tau_0,$ we have $\F^A=0$ and $\Phi^{A,B}_\cdot=0$ on $L^p([0,\al_0],\ell^1)$. Then
	\begin{align*}
	I_{\mathscr{U}}-\F^{\mathscr{A}}=\begin{pmatrix}I_{\ell^1}& -\F^Q\\ 0& I_{\ell^1}\end{pmatrix}\quad\text{on}\quad L^p([0,\al_0],\mathscr{U}).
	\end{align*}
	Thus $I_{\mathscr{U}}-\F^{\mathscr{A}}$ has a uniformly bounded inverse in $L^p([0,\al_0],\mathscr{U}),$ so $I_{\mathscr{U}}$ is an admissible feedback. Now Theorem \ref{T.1} shows that $\mathfrak{A}$ generates a strongly continuous semigroup $\mathfrak{T}:=(\mathfrak{T}(t))_{t\ge 0}$ on $\mathscr{X}$.
\end{proof}
\begin{lemma}\label{spectrum}
	Let assumptions of Theorem \ref{T.2} be satisfied. Then $\la\in \rho(\mathfrak{A})$ if and only if $1\in\rho(\A_\la),$ where we set
	\begin{align*}
	\A_\la=M D_\la+Ld_\la D_\la.
	\end{align*}
\end{lemma}
\begin{proof} According to Theorem \ref{T.1} and the proof of Theorem \ref{T.2}, $\la\in\rho(\mathfrak{A})$ if and only if $1\in\rho(\widetilde{\mathscr{M}}\mathscr{D}_\la)$. Using the notation of the proof of Theorem \ref{T.1}, we have
	\begin{align}\label{Martrice-nice}
	I_{\mathscr{U}}-\widetilde{\mathscr{M}}\mathscr{D}_\la=\begin{pmatrix} I_{\ell^1}-M D_\la & -Ld_\la\\ -D_\la& I_{\ell^1}\end{pmatrix},\quad \la\in\C.
	\end{align} This ends the proof. \end{proof}
The following is the second main result of this section.
\begin{theorem}\label{Main-section3}
	Let assumptions of Theorem \ref{T.2} be satisfied. Then the boundary control transport delay system $\mathsf{(vdp)}$ is equivalent to the following well-posed open loop system
	\begin{align}\label{open-loop}
	\begin{cases} \dot{w}(t)=\mathfrak{A}_{-1}w(t)+\mathscr{B}^K u(t),& t\ge 0,\cr w(0)=(\begin{smallmatrix}g\\\varphi\end{smallmatrix}),
	\end{cases}
	\end{align}
	where
	\begin{align*}
	\mathscr{B}^K: \C^N\to \mathscr{X}_{-1,\mathfrak{A}},\quad \mathscr{B}^K=\Big(\begin{smallmatrix} BK\\ 0\end{smallmatrix}\Big).
	\end{align*}
	Furthermore, the mild solution of is given by
	\begin{align}\label{mild}
	w(t)=\mathfrak{T}(t)(\begin{smallmatrix}g\\ \varphi\end{smallmatrix})+\Phi^{K}_t u,
	\end{align}
	for any $t\ge 0$ and $u\in L^p(\R^+,\C^N),$ where
	\begin{align*}
	\Phi^{K}_t u:=\int^t_0 \mathfrak{T}_{-1}(t-s)\mathscr{B}^K u(s)ds,\quad t\ge 0.
	\end{align*}
\end{theorem}
\begin{proof}
	We know that the boundary control transport delay system $\mathsf{(vdp)}$ is equivalent to the free-delay boundary value problem \eqref{Sy.6}. Now the result follows by combining Theorem \ref{T.2}, its proof, and Theorem \ref{nonhomo-bound}.
\end{proof}
\begin{remark}\label{autre-expression}
	Let assumptions of Theorem \ref{T.2} be satisfied. Then the mild solution of $\mathsf{(vdp)}$ is also given by
	\begin{align}\label{nice-formula}
	\begin{split}
	w(t)= \mathscr{T}(t)(\begin{smallmatrix}g\\ \varphi\end{smallmatrix})+\Phi^{\mathscr{A},\mathscr{B}}_t \mathscr{C}_\Lambda w(\cdot)+ \Phi^{\mathscr{A},\mathscr{B}}_t (\begin{smallmatrix}Ku\\ 0\end{smallmatrix})
	\end{split}
	\end{align}
	for any $t\ge 0,$  $(\begin{smallmatrix}g\\ \varphi\end{smallmatrix})\in\mathscr{X}$ and $u\in L^p(\R^+,\C^N)$, where $\mathscr{C}_\Lambda$ is Yosida extension of $\mathscr{C}$ for $\mathscr{A}$. On the other hand, if we denote by $R_\Lambda$ and $L_\Lambda$ the Yosida extension of $R$ and $L$ for $A$ and $Q,$ respectively, then $D(C_\Lambda)\times D(R_\Lambda)\subset \mathscr{C}_\Lambda$ and
	\begin{align}\label{expresso}
	\mathscr{C}_\Lambda= \begin{pmatrix} C_\Lambda& R_\Lambda\\ I& 0\end{pmatrix}\quad\text{on}\quad D(C_\Lambda)\times D(R_\Lambda).
	\end{align}
	As $\Phi^{\mathscr{A},\mathscr{B}}_t={\rm diag}(\Phi_t^{A,B},\Phi^{Q,\beta}_t)$, then by Lemma \ref{B}, Lemma \ref{regular-shift-triple}, and Theorem \ref{weiss-rep}, we deduce that $w(t)\in D(C_\Lambda)\times D(R_\Lambda)$. Now if we denote by $z(t)=w(t)_{| X}$ and $\varrho(t)=w(t)_{| L^p([-r,0],X)}$, then by combining \eqref{nice-formula} and \eqref{expresso}, we deduce that
	\begin{align}\label{z}
	\begin{split}
	z(t)&=T(t)g+\Phi^{A,B}_t [C_\Lambda z(\cdot)+R_\Lambda \varrho(\cdot)]+\Phi^{A,B}_t Ku\cr \varrho(t)&=S(t)\varphi+\Phi^{Q,\beta}_t z(\cdot).
	\end{split}
	\end{align}
	Thus $\varrho$ is the solution of an open system $(Q,\beta)$ with control function $t\mapsto z(t)$. Thus, by \cite{HI}, $\varrho$ is exactly the history function of $z(\cdot)$, so that $\varrho(t)=z_t$ for $t\ge 0$. This fact, reforce the relationship between systems $\mathsf{(vdp)}$ and \eqref{open-loop}.
\end{remark}
\section{Approximate controllability criteria}\label{S.4}
We are concerned with the approximate controllability of the transport system $\mathsf{(vdp)}$. According to Theorem \ref{Main-section3} it is preferable to study such a control property for the linear distributed system \eqref{open-loop}. Before doing so, let us first define what we means by the approximate controllability of the system \eqref{open-loop}.
\begin{definition}\label{types-controllability}Let assumptions of Theorem \ref{T.2} be satisfied. Define the control maps
	\begin{align*}\label{control-maps-closed}
	\Phi^K_t u:=\int^t_0 \mathfrak{T}_{-1}(t-s)\mathscr{B}^K u(s)ds,\quad t\ge,\; u\in L^p(\R^+,\C^N).
	\end{align*}
	Now define the space of reachable states from the initial state zero in time $t>0,$
	\begin{align*}
	\mathscr{R}^t:=\left\{\Phi^K_t u: u\in L^p([0,t],\C^N)\right\}
	\end{align*}
	The open-loop system \eqref{open-loop} is $\Sigma$-approximately controllable if and only if the set $\bigcup\{\mathscr{R}^{t}_{|\Sigma}:t>0\}$ is dense in $\Sigma$ with $\Sigma$ can be can be one of the spaces $X,L^p([-r,0],X)$ or $\mathscr{X}$.
\end{definition}
We have the following controllability result.
\begin{theorem}\label{Main-controllanility1}
	Let assumptions of Theorem \ref{T.2} be satisfied. Then the system \eqref{open-loop} (hence $\mathsf{(vdp)}$) is $\mathscr{X}-$approximately controllable if and only if for a large $\la>0$, and for any $g^\ast\in X'$ and $\varphi^\ast\in L^q([-r,0],X'),$ we have
	\begin{align*}
	& \langle D_\la (I_{\ell^1}-\A_\la)^{-1}Ku,g^\ast \rangle +
	\left\langle d_\la  D_\la (I_{\ell^1}-\A_\la)^{-1}Ku,\varphi^\ast \right\rangle=0
	\end{align*}
	for any $u\in\C^N,$ implies that $g^\ast=0$ and $\varphi^\ast=0$.
\end{theorem}
\begin{proof}
	By using the concept of closed-loop systems (see e.g. \cite{WR}), for a large $\la>0$, the Lapace transform of $t\mapsto\Phi^K_t v$ is given by
	\begin{align*}
	\widehat{\Phi^K_{\cdot} v}&= R(\la,\mathfrak{A}_{-1})\mathscr{B}^K\hat{v}(\la)\cr &=\mathscr{D}_\la (\begin{smallmatrix}K \hat{v}(\la)\\ 0\end{smallmatrix})+(I-\mathscr{C}\mathscr{D}_\la)^{-1} \mathscr{C}\mathscr{D}_\la (\begin{smallmatrix}K \hat{v}(\la)\\ 0\end{smallmatrix}).
	\end{align*}
	A similar argument as in \cite[Proposition 3]{EHR}, the system \eqref{open-loop} is $\mathscr{X}$-approximately controllable if and only if for any $u\in \C^N,$ for any $\zeta^\ast=(\begin{smallmatrix}g^\ast \\\varphi^\ast\end{smallmatrix})\in \mathscr{X}'=X'\times L^q([-r,0],X'),$
	\begin{equation*}
	\langle \mathscr{D}_\la (\begin{smallmatrix}K u\\ 0\end{smallmatrix})+(I_{\mathscr{U}}-\mathscr{C}\mathscr{D}_\la)^{-1} \mathscr{C}\mathscr{D}_\la (\begin{smallmatrix}K u\\ 0\end{smallmatrix}),\zeta^\ast \rangle=0
	\end{equation*}
	implies that $\zeta^\ast=0$, so that $g^\ast=0$ and $\varphi^\ast=0$. The result follows by computing the inverse of the operator $I_{\mathscr{U}}-\mathscr{C}\mathscr{D}_\la$ given by \eqref{Martrice-nice}.
\end{proof}
\begin{proposition}\label{partial-prop} Let assumptions of Theorem \ref{T.2} be satisfied. Then we have
	\begin{itemize}
		\item[{\rm(i)}] The system \eqref{open-loop} (or $\mathsf{(vdp)}$) is $X$-approximately controllable if and only if for any $u\in \C^N$ and $g^\ast\in X'$, \begin{align}\label{F2}
		\langle  D_\la (I_{\ell^1}-\A_\la)^{-1}Ku,g^\ast \rangle=0
		\end{align}
		implies that $g^\ast=0$.
		\item [{\rm(ii)}] The system \eqref{open-loop} (or $\mathsf{(vdp)}$) is $L^p([-r,0],X)$-approximately controllable if and only if for any $u\in \C^N$ and $\varphi^\ast\in L^q([-r,0],X')$, \begin{align}\label{F3}
		\left\langle d_\la  D_\la (I_{\ell^1}-\A_\la)^{-1}Ku,\varphi^\ast \right\rangle=0
		\end{align}
		implies that $\varphi^\ast=0$.
	\end{itemize}
\end{proposition}
\begin{proof}
	Let $w(t,0,u)$ be the solution of  \eqref{open-loop} for the initial conditions $g=0$ and $\varphi=0$, this  is exactly the control maps $\Phi^K_tu$. Then the Laplace transform of the projection of    $\Phi^K_tu$ is exactly the Laplace transform of the function $t\mapsto z(t,0,u)$, where $z(\cdot)$ is given by \eqref{z}. Taking Laplace form in both sides of this formula, for large $\la,$ and using \cite{HIR}, one can see that
	\begin{align*}
	\widehat{z(\cdot,0,u)}(\la) &=(I_{\ell^1}-D_\la(R_\Lambda+Ld_\la))^{-1}D_\la Ku\cr &=D_\la Ku+D_\la (I_{\ell^1}-\A_\la)^{-1}(R D_\la+Ld_\la D_\la)K\hat{u}(\la)\cr &= D_\la (I+(I_{\ell^1}-\A_\la)^{-1} \A_\la K\hat{u}(\la)
	\cr &= D_\la(I_{\ell^1}-\A_\la)^{-1}K\hat{u}(\la)
	\end{align*}
	On other hand,  Laplace transform of the projection of    $\Phi^K_tu$ is exactly the Laplace transform of the function $t\mapsto z_t(\cdot,0,u)=\Phi^{Q,\beta}z(\cdot),$ see \eqref{z}. Thus $\widehat{z_{\cdot}(0,u)}(\la)=d_\la \widehat{z(\cdot,0,u)}(\la)$. Thus the results follows from the proof of \cite[Proposition 3]{EHR}.
\end{proof}
\begin{remark}\label{relation-between-controllability}
	\begin{enumerate}
		\item Let assumptions of Theorem \ref{T.2} be satisfied. Assume that the system \eqref{open-loop} is $\mathscr{X}$-approximately controllable. By taking $g^\ast\in X'$ and $\varphi^\ast=0,$ in the formula of Theorem \ref{Main-controllanility1} and using the point (i) of Proposition \ref{partial-prop}, we deduce that \eqref{open-loop} is $X$-approximately controllable. Moreover, if we take $g^\ast=0$ and $\varphi^\ast\in L^q([-r,0],X'),$ in in the formula of Theorem \ref{Main-controllanility1} and using the point (i) of Proposition \ref{partial-prop}, we deduce that \eqref{open-loop} is $L^p([-r,0],X)$-approximately controllable.
		
		\item Let $\la>0$ sufficiently large. For simplicity we put
		\begin{align*}
		\mathbb{X}_\la=  D_\la (I_{\ell^1}-\A_\la)^{-1}K.
		\end{align*}
		Assume that the system \eqref{open-loop} is $L^p([-r,0],X)$-approximately controllable, and let $g^\ast\in X'$ et $u\in\C^N$ such that $\langle \mathbb{X}_\la u,g^\ast\rangle=0$. We select $\varphi^\ast(\theta)=e^{-\la \theta}g^\ast$ for $\theta \in [-r,0]$. Thus $\varphi^\ast\in L^q([-r,0],X')$. We have
		\begin{align*}
		\langle d_\la \mathbb{X}_\la u,\varphi^\ast\rangle= -r \langle \mathbb{X}_\la u,g^\ast\rangle=0.
		\end{align*}
		Then by Proposition \ref{relation-between-controllability} we have $\varphi^\ast=0$, so that $g^\ast=0$. This means that the condition \eqref{F3} implies \eqref{F2}. We conclude that $L^p([-r,0],X)$-approximately controllable of the system \eqref{open-loop} implies the $X$-approximately controllable of the system \eqref{open-loop}.
	\end{enumerate}
\end{remark}

The following theorem gives a practical characterization for approximate controllability of vertex delay problems $\mathsf{(vdp)}$.
\begin{theorem}\label{T.4}
	Let assumptions of Theorem \ref{T.2} be satisfied and assume that $\inf_{j\in\N}\tau_j(0,1)>0$. There is $\mu_0>0$ such that for any ${\rm Re}\la> \mu_0$, the vertex delay problem $\mathsf{(vdp)}$ is $X$-approximately controllable if and only if
	\begin{align}\label{4.5}
	\text{Cl}\left(\text{span } \left\{\mathbb{A}_{\la}^{k}K,\; k=0,1,2,\cdots \right\}\right)=\ell^{1},
	\end{align}
	where, $ \text{Cl}(\cdot) $ means the closure of a set.
\end{theorem}
\begin{proof}
	Let us first show that $\|\A_\la\|\to 0$ as ${\rm Re}\la\to +\infty$. Inspecting from  \cite[Lemma 6.1]{Hadd-SF}, we have  $\|Rd_\la\|\to 0$ as ${\rm Re}\la\to+\infty$. On the other hand, from \cite[chap.3]{TW}, as $(A,B)$ is admissible then for any $\omega>\om_0(A)$ (semigroup type) \begin{align*}\|D_\la\|\le \frac{b}{({\rm Re}\la-\om)^{\frac{1}{q}}} \end{align*}
	for any ${\rm Re}\la>\omega$ and for a constant $b>0$. Thus $\|Rd_\la D_\la\|\to 0$ as ${\rm Re}\la\to+\infty$. From \eqref{Dirichlet-D}, we have
	\begin{align*}
	MD_\la=\Upsilon \left({\rm diag}(e^{\xi_j(0,1)-\la \tau_j(0,1)})\right)_{j\in\N}.
	\end{align*}
	For any $j\in \N,$ and ${\rm Re}\la>0$ we have $e^{-{\rm Re}\la \tau_j(0,1)}\le e^{-{\rm Re}\la\inf_{j\in\N}\tau_j(0,1)}\to 0$ as ${\rm Re}\la\to +\infty,$ due to $\inf_{j\in\N}\tau_j(0,1)>0$. This implies that $\|MD_\la\|\to 0$ as ${\rm Re}\la\to 0$, then we proved the first claim. Then there exists a large $\mu_0>0$ such that $ \Vert \mathbb{A}_\la\Vert <1 $ for any ${\rm Re}\la>\mu_0$. Thus, using the Neumann series,
	\begin{align*}
	D_{\la}(I_{\ell^{1}}-\mathbb{A}_{\la})^{-1}K=D_{\la}\sum_{k\geqslant 0}\mathbb{A}_{\la}^{k}K,\qquad {\rm Re}\la>\mu_0.
	\end{align*}
	According to Proposition \ref{partial-prop}  and \eqref{Dirichlet-D}, the vertex delay problem $\mathsf{(vdp)}$ is $X$-approximately controllable if and only if
	$$ \text{Rg}\left(\text{diag}(e^{\xi_{j}(\cdot ,1)-\la\tau_{j}(\cdot ,1)})_{j\in \mathbb{N}}\sum_{k\ge 0}\mathbb{A}_{\la}^{k}K\right) $$
	is dense in $ X $. On the other hand, it is well known
	that
	\begin{align*}
	X:=L^{p}([0,1],\ell^{1})=Cl\big(L^{p}([0,1],\mathbb{C})\otimes \ell^{1}\big).
	\end{align*}
	As a matter of fact, the vertex delay problem $\mathsf{(vdp)}$ is $X$-approximately controllable if and only if the condition \eqref{4.5} fulfill, since by the Stone-Weierstrass theorem, for any $j\in \mathbb{N}$,
	$$ \text{span }\left\{ e^{-\la\tau_{j}(\cdot ,1)}:\; \Re e\, \la> \mu_0 \right\}$$
	is dense in $ C([0,1],\mathbb{C})$, and hence
	$$ \text{span }\left\{ e^{\xi_{j}(\cdot ,1)-\la\tau_{j}(\cdot ,1)}:\; \Re e\, \la> \mu_0 \right\}$$ is also dense in
	$  C([0,1],\mathbb{C})\subset L^{p}([0,1],\mathbb{C}) $. This complete the proof.
\end{proof}
The following result  generalizes the corresponding on in \cite[Theorem 4.5]{EHR} from finite to  infinite network.
\begin{corollary}\label{C4.Cor2}
	Let the assumptions {\bf  (A1)} to {\bf  (A3)} be satisfied and assume that $\sup_{i\in\N}\|c_i\|_\infty<\infty$ and  $\t_0:=\inf_{ j\in \mathbb{N}}\t_j(0,1)>0$. Then, there is $\mu_0>0$ such that for any ${\rm Re} \la> \mu_0$, the transport network system $\mathsf{(vdp)}$ without delay, i.e. $L=0$, is approximately controllable if and only if
	\begin{align}\label{4.9}
	\text{Cl}\left(\text{span} \left\{\mathbb{A}_{\la}^{k}K,\; k=0,1,2,\cdots \right\}\right)=\ell^{1},
	\end{align}
	where
	\begin{align*}
	\mathbb{A}_{\la}=MD_{\la}=c(1)^{-1}\mathbb{B}c(0)\;{\rm diag}(e^{\xi^{j}(0,1)-\la \tau^{j}(0,1)})_{j\in \N}.
	\end{align*}
\end{corollary}
\begin{remark}
	Observe that the algebraic condition \eqref{4.9} in Corollary \ref{C4.Cor2} does not require the condition $q_j\le 0$, compare with the corresponding result is \cite[Theorem 4.2]{EHR}.
\end{remark}

\section{Approximate controllability of an air traffic flow management}\label{S.5}
In this section, we study a particular case of the delay network equation $\mathsf{(vdp)}$. In fact, the problem of regulating air traffic flows in a region of airspace (airway) is considered. This problem is known as the \emph{Air Traffic Flow Management} (ATFM), see e.g. \cite{SSB}.

Consider a finite connected graph $ \mathsf{G} $ having $ v_{1},\ldots, v_{n}$ vertices and $ e_{1},\ldots ,e_{m}$ edges and let a flow of air traffic on this graph. According e.g. to \cite{MSB}, a such traffic is described by the following system of transport equations
\begin{align}\label{S4.1}
\begin{cases}
\dfrac{\partial }{\partial t}z_{j}(t,x)= c_{j}(x)\dfrac{\partial }{\partial x}z_{j}(t,x),& x\in (0,1),\; t\geq 0, \\
z_{j}(0,x)= g_{j}(x),& x\in (0,1), \\
c_{j}(1)z_{j}(t,1)= \sum_{k=1}^{m} h_{jk}\big[c_{k}(0)z_{k}(t,0)+\int_{0}^{1}c_k(x)z_k(t-r,x)dx\big]+k_{jl}u_l(t),& t\geq 0,\\
z_{j}(\theta,x)=\varphi_{j}(\theta,x),& \theta\in[-r,0], x\in (0,1),
\end{cases}
\end{align}
for $ j=1,\ldots, m $ and $l=1,,\ldots, n_0$ with $n_0\leq n$, where $c_j\in L^\infty([0,1])$ are such that $c_j(x)\ge c_1$ for all $x\in[0,1]$, $1\leq j\leq m$ and a constant $c_1>0$, . The coefficients $ h_{jk} $ denotes the entries of the so-called (transposed) \emph{allocation matrix} $ \mathbb{H}  $, where $ 0\leq h_{jk}\leq 1 $ determines the proportion of aircrafts arriving from edge $ e_k $ leaving into edge $ e_j $. If $ \B=(b_{jk}) $ denotes the (transposed) adjacency matrix of the line graph, we set $ h_{jk}=0 $ if $ b_{jk}=0 $ and assume that
\begin{align}\label{S4.2}
\sum_{j=1}^{m} h_{jk}=1, \; \; \; \forall k=1,\ldots,m.
\end{align}
\begin{remark}\label{RF}
	It should be noted that in ATFM the traffic is usually considered
	on an isolated junction of the network, and it is modeled by a system of transport equations, where the edges are identified with $ n+m $ intervals (i.e., $ m $-incoming edges and $ n $-outgoing edges), with connection described by the junction allocation matrix $ \mathbb{H}$, see \cite{MSB,BRT} for more details.
\end{remark}

If we denote by $K:=(k_{jl})$, then Theorem \ref{T.4} shows that the linear Eulerian model \eqref{S4.1} is boundary approximately controllable if and only if, for $ {\rm Re} \mu>0 $,
\begin{align}\label{S4.3}
\text{Rank} \begin{pmatrix}
K &
\mathbb{A}_{\mu}K
&\ldots &
\mathbb{A}_{\mu}^{m-1}K
\end{pmatrix}=m,
\end{align}
where
\begin{align*}
\mathbb{A}_{\mu}&:=c^{-1}(1)\mathscr{H}c(0)\big[\text{diag}\left(\begin{smallmatrix}e^{-\mu \tau_{j}(0,1)}\end{smallmatrix}\right)_{j=1,\ldots, m}
+c^{-1}(0)\text{diag}\left(\begin{smallmatrix}\int_{0}^{1}e^{-\mu (r+\tau_{j}(x,1))}c_{j}(x)dx\end{smallmatrix}\right)_{j=1,\ldots, m} \big].
\end{align*}
\begin{remark}\label{RF1}
	Note that the introduced model \eqref{S4.1} is a variant of the transport network system considered in \cite[Section 5]{BDK}, where the transposed matrix $ \mathbb{H}^{\top} $ of every junction allocation matrix $ \mathbb{H} $ is a submatrix of the bigger matrix $ \mathscr{H} $, which, according to the linear Eulerian model \eqref{S4.1}, contain full information on the aggregate aircraft flow of each airspace of the entire network. Further, we consider the situation of airborne delay.
\end{remark}
As noted in Remark \ref{RF}, the air traffic flows are studied in an isolated junction of the network, which may cause a lack of some numerical values of the matrix $ \mathscr{H} $. In such cases, it is difficult to verify the Kalman-type condition  \eqref{S4.3}. However, based on the concept of structural controllability \cite{SP}, the following result provide a useful controllability criteria for the controllability of boundary controlled linear Eulerian models.
\begin{theorem}
	System \eqref{S4.1} is approximately controllable if and only if, for any $\mu\in\C$ such that $ {\rm Re}\mu>0 $, the following matrix
	$$
	\mathfrak{C}:=\left(\begin{smallmatrix}
	K & I & 0 & 0 & 0 & \cdots & \cdots & & &  0\\
	0 & -\mathbb{A}_\mu & K & I & 0 & \cdots & \cdots& & &  0\\
	0 & 0 & 0 & -\mathbb{A}_\mu & K & \cdots & \cdots &  &  &\\
	\vdots & \vdots & \vdots & \vdots & \vdots & \vdots & \vdots  &  &  &\vdots\\
	0 & 0 & 0 & & \cdots & & -\mathbb{A}_\mu & K & I & 0 \\
	0 & 0 & 0 & & \cdots & & 0 & 0 & -\mathbb{A}_\mu & K
	\end{smallmatrix}\right),
	$$
	is not of form $ (n^{2}) $.
\end{theorem}
\begin{proof}
	By virtue of Theorem \ref{T.4} the linear Eulerian model \eqref{S4.1} is approximately controllable if and only if \eqref{S4.3} holds for any $\mu\in\C$ such as $ {\rm Re}\mu>0 $. Equivalently, \eqref{S4.1} is approximately controllable if and only if the extended controllability matrix $ \mathfrak{C} $ has rank $ n^{2} $.  \cite[Theorem 4.1]{SP} yields the result.	
\end{proof}

\begin{remark}
	Notice that the above result ensures the existence of a feasible solution for the linear Eulerian model of ATFM problem, i.e., in the sense of finding the best combination of flow controls and departure rates for the controllability of a given air traffic environment \cite[Sec. II.B.1]{MSB}. Moreover the above controllability criteria is more robust and of practical interest for linear Eulerian models. This is because it is independent of the detailed values of the aggregate aircraft flow for each airspace in the air traffic environment.
\end{remark}
\section{Conclusion} In this paper, we have presented a semigroup approach to the well-posedness and approximate controllability of a transport tree-like network with infinitely many edges including hereditary effects in the transmission conditions. Based on the transformation of the delay system into a delay-free open loop system, we first introduced three concepts of controllability and showed the relationship  between these concepts. Furthermore, by using a Laplace transform technique, we characterised each of those controllability.  In particular, we have a established a controllability criterion in terms of a Kalman-type rank condition involving the graph structure. The approximate controllability of an air traffic flow management with a  discrete delay at the boundary conditions is also considered.


\section{Appendix}
Consider the following linear time-invariant system
\begin{align*}
\Sigma(\mathbb{A},\mathbb{K}):\; \; \; \; \dot{x}(t)=\mathbb{A} x(t)+\mathbb{K} u(t),\;\; t\geq 0,
\end{align*}
with $ x(t)\in \mathbb{R}^{n} $ and $ u(t)\in \mathbb{R}^{m} $, and such as $ \mathbb{A} $ and $ \mathbb{K} $ are \emph{structured} matrices in the sense that the elements of the above matrices are \emph{zero/nonzero} type.
\begin{definition}
	The two systems $ \Sigma(\mathbb{A},\mathbb{K}) $ and $ \Sigma'(\mathbb{A}',\mathbb{K}') $ are structurally equivalent if there is a one-to-one correspondence between the locations of their fixed zero and nonzero entries.
	
	A system $\Sigma(\mathbb{A},\mathbb{K}) $ is called structurally controllable if there exists a system structurally equivalent to $ \Sigma(\mathbb{A},\mathbb{K}) $ which is controllable in the usual sense, i.e., fulfilled Kalman's controllability rank condition.
\end{definition}
\begin{example}\label{A.1}
	Each of the following matrices is structural
	\begin{align*}
	Q_{0}:=\left(\begin{smallmatrix}
	0 & 0 & 0 & 0 & 0\\
	0 & 0 & 0 & 0 & 0\\
	x & x & x & x & x\\
	x & x & x & x & x\\
	x & x & x & x & x
	\end{smallmatrix}\right),\; Q_{1}:=\left(\begin{smallmatrix}
	x & 0 & 0 & 0 & 0\\
	x & 0 & 0 & 0 & 0\\
	x & 0 & 0 & 0 & 0\\
	x & x & x & x & x\\
	x & x & x & x & x
	\end{smallmatrix}\right).
	\end{align*}
\end{example}

Note that if a system is structurally controllable then it is controllable for almost all parameter values except for those lie on a proper algebraic variety in the parameter space. This shows that structural controllability is a generic property of the system, since a proper algebraic variety has Lebesgue measure zero \cite[Proposition 3.1]{SP}.

We now introduce the basic technical condition for studying structural controllability.
\begin{definition}
	An $ (n\times s) $ matrix $A\; (s\geq n )$ is said to be of form $ (t) $ for some $ 1 \leq t \leq n $ if, for some $ k $ in the range $ s-t< k \leq s $, $ A$ contains a zero sub-matrix of order $ (n+s-t-k+1)\times k$.
\end{definition}
\begin{example}
	Both matrices from Example \ref{A.1} are of form $ (4) $, but $ k=5 $ for the matrix $ Q_{0} $ while $ k=4 $ for the matrix $ Q_{1} $.
\end{example}
\begin{lemma}
	For any $ t $, $ 1 \leq t \leq n $, $ \text{Rank } A<t $ for every $ p\in \mathbb{R}^{n} $ if and only if $ A $ has form $ (t) $.
\end{lemma}

The materiel presented in this appendix is taken from Shields and Pearso \cite{SP}, where they extended Lin's result \cite{Li} (in the case of single-input) on the structural controllability of a given structured $ \Sigma(\mathbb{A},\mathbb{K}) $ system to the multi-input case.

\end{document}